\documentclass[journal, onecolumn,letterpaper,draftcls, 12 pt]{IEEEtran} 
\usepackage{times} 
\usepackage{amsmath,lipsum} 
\usepackage{amssymb,amsfonts}  
\usepackage{latexsym,theorem,epsfig}
\usepackage{graphicx} 
\usepackage{dsfont}
\usepackage{mathtools}
\usepackage{subfigure}
\usepackage{lscape}
\usepackage{color}
\usepackage{float}
\usepackage{enumerate}
\usepackage[usenames,dvipsnames]{pstricks}
\usepackage{epsfig}
\usepackage{pst-grad} 
\usepackage{pst-plot} 
\usepackage{mathrsfs}
\usepackage{extarrows}
\usepackage{caption}
\usepackage[percent]{overpic}

\newtheorem{theorem}{Theorem}

\newtheorem{assumption}{Assumption}

\newtheorem{definition}{Definition}
{\theorembodyfont{\rmfamily} }
\newtheorem{lemma}{Lemma}

\newtheorem{proposition}{Proposition}

\allowdisplaybreaks
\IEEEoverridecommandlockouts

\DeclareMathOperator*{\argmin}{arg\,min}
\title{Distributed Learning with Infinitely Many Hypotheses}

\author{Angelia Nedi\'{c}, Alex Olshevsky and C\'{e}sar A.\ Uribe
\thanks{The authors are with the Coordinated Science Laboratory, University of Illinois, 1308 West Main Street, Urbana, IL 61801, USA, \{angelia,aolshev2,cauribe2\}@illinois.edu.
	This research is supported partially by the National Science Foundation under
	grants no.\ CCF 11-11342 and no.\ CMMI-1463262 and by the Office of Naval Research under grant
	no.\ N00014-12-1-0998.}  
}
\begin{document}
\maketitle
\begin{abstract}
We consider a distributed learning setup where a network of agents sequentially access realizations of a set of random variables with unknown distributions. The network objective is to find a parametrized distribution that best describes their joint observations in the sense of the Kullback-Leibler divergence. Apart from recent efforts in the literature, we analyze the case
of countably many hypotheses and the case of a continuum of hypotheses. 
We provide non-asymptotic bounds for the concentration rate of the agents' beliefs around the correct hypothesis in terms of the number of agents, the network parameters, and the learning abilities of the agents. 
Additionally, we provide a novel motivation for a general set of distributed Non-Bayesian update rules as instances of the distributed stochastic mirror descent algorithm. 
\end{abstract}

\section{Introduction}

Sensor networks have attracted massive attention in past years due to its extended range of applications and its ability to handle distributed sensing and processing for systems with inherently distributed sources of information, 
e.g., power networks, social, ecological and economic systems, surveillance, disaster management health monitoring, etc.\   \cite{gun10,mai02,cho03}. 
{For such distributed systems, one can assume complete communication between every source of information (e.g. nodes or local processing unit) and centralized processor can be cumbersome.} 
Therefore, one might consider cooperation strategies where nodes with limited sensing capabilities distributively aggregate information to perform certain global estimation or processing task.  

Following the seminal work of Jadbabaie et al.\ in \cite{jad12,mol15}, there have been many studies of Non-Bayesian rules for distributed algorithms. 
Non-Bayesian algorithms involve an aggregation step, usually consisting of weighted geometric or arithmetic average of the received beliefs, and a Bayesian update that is based on the { locally available data.} Therefore, one can exploit results from consensus literature~\cite{ace08,tsi84,jad03,ned13,ols14} and Bayesian learning 
literature~\cite{ace11,mos14}. 
Recent studies have proposed several variations of the Non-Bayesian approach and have proved consistent, geometric and non-asymptotic convergence rates for a 
general class of distributed algorithms; from asymptotic analysis \cite{sha13,lal14,qip11,qip15,sha15,rah15} to non-asymptotic bounds \cite{sha14,ned15,lal14b}, time-varying directed graphs \cite{ned15b} and transmission and node failures ~\cite{su16}. 
 
In contrast with the existing results that assume a finite hypothesis set, 
in this paper, we are extending the framework to the cases of a countable many and a {continuum of hypotheses.} 
We build upon the work in \cite{bir15} on non-asymptotic behaviors of Bayesian estimators to construct non-asymptotic concentration results for distributed learning. In the distributed case, the observations will be scattered among a set of nodes or agents and the learning algorithm should guarantee that every node in the network will learn the correct parameter as if it had access to the complete data set. 
Our results show that in general the network structure will induce a transient time after which all agents will learn at a network independent rate, where the rate is geometric.

{\it The contributions of this paper are threefold}: First, we provide an interpretation of a general class of distributed Non-Bayesian algorithms as specific instances of a distributed version of the stochastic mirror descent. 
This motivates the proposed update rules and makes a connection between the Non-Bayesian learning literature in social networks and the Stochastic Approximations literature. Second, 
we establish a non-asymptotic concentration result for the proposed learning algorithm when the set of hypothesis is countably infinite. 
Finally, we provide a non-asymptotic bound for the algorithm when the hypothesis set is a bounded subset of 
$\mathbb{R}^d$. This is an initial approach to the analysis of distributed 
Non-Bayesian algorithms for a more general family of hypothesis sets.

This paper is organized as follows: Section~\ref{setup} describes the studied problem and the proposed algorithm, together
with the motivation behind the proposed update rule and its connections with distributed stochastic mirror descent algorithm. Section~\ref{sec_count} and Section~\ref{sec_cont} provide the non-asymptotic concentration rate results for 
the beliefs around the correct hypothesis set 
for the cases of countably many and continuum of hypotheses, respectively. 
Finally, conclusions are presented in Section~\ref{sec_conclusion}. 
 
 \noindent
{\bf Notation}:
The set $B^c$ denotes the complement of a set $B$. 
Notation $\mathbb{P}_B$ and $\mathbb{E}_B$ denotes 
the probability measure and expectation under a distribution $P_B$. 
The $ij$-th entry of a matrix $A$ 
is denoted by  $\left[A\right]_{ij}$ or $a_{ij}$. 
Random variables are denoted with upper-case letters, while the corresponding lower-case letters denote their realizations. 
Time indices are indicated by subscripts and the letter $k$. 
Superscripts represent the agent indices, which are usually $i$ or~$j$. 
 
 \section{Problem Formulation }\label{setup}
 
We consider the problem of distributed non-Bayesian learning, where a network of agents access 
sequences of realizations of a random variable with an unknown distribution.
The random variable is assumed to be of finite dimension with the constraint that each agent can access only a strict subset of the entries of the realizations (e.g., an $n$-dimensional vector and $n$ agents each observing a single entry). 
Observations are assumed to be independent among the agents. 
We are interested in situations where no single agent has the ability to learn the underlying distribution from
its own observations, while collectively the agents can do so if they collaborate.
The learning objective is for the agents to jointly agree on a distribution 
(from a parametrized family of distributions or a hypothesis set) 
that best describes the observations in a specific sense (e.g., Kullback-Leibler divergence).
Therefore, the distributed learning objective requires collaboration among the agents which 
can be ensured by using some protocols for information aggregation and coordination. 
Specifically in our case, agent coordination consists of sharing their estimates  (\textit{beliefs}) of 
the best probability distribution over the hypothesis set.
 
\begin{figure}[ht]
	 	\centering
	 	\captionsetup{justification=centering}
	 	\includegraphics[width=0.4\textwidth]{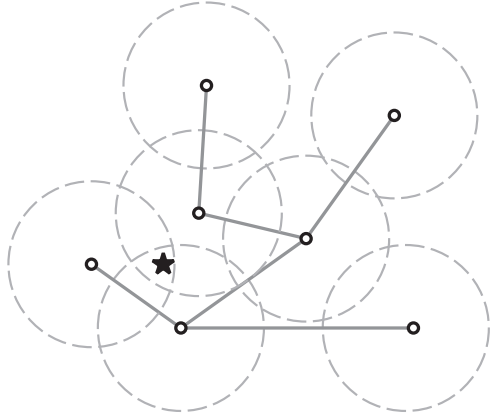}
	 	\caption{Distributed source localization example}
	 	\label{location}
	 \end{figure}

Consider, for example, the distributed source location problem with limited sensing capabilities \cite{rab04,rab05}. In this scenario a network of $n$ agents receives noisy measurements of the distance to a source, 
where sensing capabilities of each sensor might be limited to a certain region. 
The group objective is to jointly identify the location of the source and that every node knows the source location. 
Figure \ref{location} shows an example, where a group of $7$ agents (circles) wants to localize a source (star). 
There is an underlying graph that indicates the communication abilities among the nodes. 
Moreover, each node has a sensing region indicated by the dashed line around it. 
Each agent $i$ obtains realizations of the random variable $S_k^i = \|x^i - \theta^*\| + W_k^i$, 
where $\theta^*$ is the location of the source, $x^i$ is the position of agent $i$ and 
$W_k^i$ is a noise in the observations.
If we consider $\Theta$ as the set of all possible locations of the source, then each $\theta \in \Theta$ will induce a probability distribution about the observations of each agent. 
Therefore, agents need to cooperate and share information in order to guarantee that all of them correctly localize the target.

We will consider a more general learning problem,
where agent observations are drawn from an unknown joint distribution $\boldsymbol{f} = \prod_{i=1}^{n} f^i$, 
where $f^i$ is the distribution governing the observations 
of agent $i$. 
We assume that $\boldsymbol{f}$ is an element of $\mathscr{P} = \prod_{i=1}^n \mathscr{P}^i$, the space of all joint probability measures 
for a set of {$n$ independent random variables} $\{S^i\}_{i=1}^n$ (i.e., $S^i$ is distributed according to an unknown distribution $f^i$). Also, we assume that each $S^i$ takes values in a finite set. When these random variables are considered at time $k$, we denote them by $S^i_k$.

	Later on, for the case of countably many hypotheses, we will use 
	the {\it pre-metric} space $\left(\mathscr{P},D_{KL}\right)$ as the vector space $\mathscr{P}$ equipped with 
	the Kullback-Liebler divergence. This will generate a topology, where we can define an open ball 
	$\mathcal{B}_r(p)$ with a radius $r>0$ centered at a point $p \in \mathscr{P}$ by 
	${\mathcal{B}_r(p) = \{q \in \mathscr{P}| D_{KL}(q,p) < r \}}$. 
	When the set of hypothesis is continuous, we instead equip $\mathscr{P}$ with the Hellinger distance~$h$
	to obtain the {\it metric} space $\left(\mathscr{P},h\right)$, 
	which we use to construct a special covering of subsets $B\subset \mathscr{P}$ consisting of $\delta$-separated sets.

Each agent constructs a set of hypothesis parametrized by $\theta \in \Theta$ about the distribution $f^i$. Let $\mathscr{L}^i = \{P^i_\theta | \theta\in \Theta\}$ be a parametrized family of probability measures for $S_k^i$ 
with densities $\ell^i\left(\cdot|\theta\right) = dP^i_\theta / d\lambda^i$ with respect to a dominating measure\footnote{A measure $\mu$ is dominated by $\lambda$ if $\lambda(B) = 0$ implies $\mu(B)=0$ for every measurable set $B$.}  
$\lambda^i$. Therefore, the learning goal is to distributively solve the following  problem:
 \begin{align}\label{opt_problem}
 \min_{\theta \in \Theta} F(\theta) & \triangleq D_{KL}\left(\boldsymbol{f}\|\boldsymbol{\ell}\left(\cdot | \theta\right)\right) \ \ \   \\
 & = \sum\limits_{i=1}^n D_{KL}\left(f^i\|\ell^i\left(\cdot|\theta\right)\right) \nonumber
 \end{align}
 where $\boldsymbol{\ell}\left(\cdot | \theta\right)=\prod_{i=1}^n \ell^i\left(\cdot|\theta\right)$ and 
 $D_{KL}\left(f^i\|\ell^i\left(\cdot|\theta\right)\right)$ is 
 the Kullback-Leibler divergence between the true distribution of $S_k^i$ and $\ell^i( \cdot | \theta)$ 
 that would have been seen by agents $i$ if hypothesis $\theta$ were correct. 
For simplicity  {\it we will assume that there exists a single $\theta^* \in \Theta$} 
such that $\ell^i\left(\cdot|\theta^*\right) = f^i$ almost everywhere for all agents. 
Results readily extends to the case when the assumption does not hold 
(see, for example,~\cite{ned14,ned15,ned15b} which disregard this assumption).

The problem in Eq.~\eqref{opt_problem} consists of finding the parameter $\theta^*$ such that 
$\boldsymbol{\ell}\left(\cdot | \theta\right) = \prod_{i=1}^{n}\ell^i\left(\cdot|\theta\right)$ minimizes 
its Kullback-Liebler divergence to $\boldsymbol{f}$. However, $\mathscr{L}^i$ is only available to agent $i$ and the distribution $\boldsymbol{f}$ is unknown. Agent $i$ gains information on $f^i$ by observing realizations 
$s_{k}^i$ of $S_k^i$ at every time step $k$. 
The agent uses these observations to construct a sequence $\{\mu_k^i\}$
of probability distributions over the parameter space $\Theta$. We refer to  these distributions as
agent $i$ \textit{beliefs},
where $\mu_k^i\left(B\right)$ denotes the belief, at time $k$,  that agent $i$ has about 
the event $\theta^* \in B \subseteq \Theta$  for a measurable set $B$.
 
We make use of the following  assumption. 
  \begin{assumption}\label{assum:init}
  	For all agents $i=1,\ldots,n$ we have:
  	\begin{enumerate}[(a)]
  		\item
  		There is a unique hypothesis $\theta^*$ such that $\ell^i\left(\cdot|\theta^*\right) = f^i$.
  		\item 
  		If $f^i\left(s^i\right)>0$, then there exists an $\alpha >0$  such that $\ell^i\left(s^i | \theta \right)  > \alpha$ for all $\theta \in \Theta$.
  	\end{enumerate}
  \end{assumption}
 {Assumption~\ref{assum:init}(a) guarantees that we are working on the realizable case and 
there are no conflicting models among the agents, see~\cite{ned14,ned15,ned15b}
for ways of how to remove this assumption. }
{Moreover in Assumption~\ref{assum:init}(b), the lower bound $\alpha$ assumes the set of hypothesis are dominated by $f^i$ (i.e., our hypothesis set is absolutely continuous with respect to the true distribution of the data) and provide a way to show bounded differences when applying the concentration inequality results.}

Agents are connected in a network 
$\mathcal{G} = \{V,E\}$ where $V = \{1,2,\ldots,n\}$ is the set of agents and 
$E$ is a set of undirected edges, where $\left(i,j\right) \in E$ if agents $i$ and $j$ can communicate with each other. 
If two agents are connected they share their beliefs over the hypothesis set at every time instant $k$. We will propose a distributed protocol to define how the agents update their beliefs based on their local observations and the beliefs received from their neighbors. 
Additionally, each agent weights its own belief and the beliefs of its neighbors; 
we will use $a_{ij}$ to denote the weight that agent $i$ assigns to beliefs coming from its neighbor $j$, and $a_{ii}$ 
to denote the weight that the agent assigns to its own beliefs. 
The assumption of static undirected links in the network is made for simplicity of the exposition.
The extensions of the proposed protocol to more general cases of time varying undirected and directed graphs
can be done similar to the work in \cite{ned14,ned15,ned15b}.
 
Next we present the set of assumptions on the network over which the agents are interacting.
 \begin{assumption}\label{assum:graph}
 	The graph $\mathcal{G}$ and matrix $A$ are such that:
 	\begin{enumerate}[(a)]
 		\item $A$ is doubly-stochastic with $\left[A\right]_{ij} = a_{ij} > 0$ if $(i,j)\in E$.
 		\item If $(i,j) \notin E$ for some $i \neq j$ then $a_{ij}=0$.
 		\item $A$ has positive diagonal entries, $a_{ii}>0$ for all $i \in V $.
 		\item If $a_{ij}>0$, then $a_{ij} \ge \eta$ for some positive constant $\eta$.
 		\item The graph $\mathcal{G}$ is  connected.
 	\end{enumerate}
 \end{assumption} 
 Assumption~\ref{assum:graph} is common in distributed optimization and consensus literature. 
It guarantees convergence of the associated Markov Chain and defines bounds 
on relevant eigenvalues in terms of the number $n$ of agents.
To construct a set of weights satisfying Assumptions \ref{assum:graph}, for example, one can 
consider a lazy metropolis ({stochastic}) matrix of the form 
$\bar{A} = \frac{1}{2}I + \frac{1}{2}A$, where $I$ is the identity matrix and $
A$ is a stochastic matrix whose off-diagonal entries satisfy
		\begin{align*}
		a_{ij} = \left\{  \begin{array}{l l}
		\frac{1}{\max \left\{ d^i+1,d^j+1 \right\}  } & \quad \text{if $(i,j) \in E$ }\\
		0 & \quad \text{if $(i,j) \notin E$}
		\end{array} \right. 
		\end{align*}
		where $d^i$ is the degree (the number of neighbors) of node $i$.
Generalizations of Assumption~\ref{assum:graph} 
to time-varying undirected are readily available for weighted averaging
and push-sum approaches~\cite{ned09,ols14,ned13}.	

\section{Distributed Learning Algorithm }\label{sec_stoc}
In this section, we present the proposed learning algorithm
and a novel connection between Bayesian update and the stochastic mirror descent method.
We propose the following (theoretical) 
algorithm, where each node updates its beliefs 
on a measurable subset $B\subseteq\Theta$ according to the following update rule:
for all agents $i$ and all $k\ge1$,
\begin{align}\label{protocol_cont}
\mu_{k}^i\left(B\right) & = \frac{1}{Z_{k}^i}\int\limits_{\theta \in B} \prod\limits_{j=1}^{n} \left( \frac{d\mu_{k-1}^j\left(\theta\right)}{d\lambda(\theta)}\right)^{a_{ij}}  \ell^i(s_{k}^i|\theta) d\lambda\left(\theta\right)
\end{align}
where $ Z_{k}^i $ is a normalizing constant and $\lambda$ is a probability distribution on 
$\Theta$ with respect to which every $\mu_k^j$ is absolutely continuous. 
The term  $d \mu^j_k(\theta)/d \lambda(\theta)$ is the
Radon-Nikodym derivative of the probability distribution $\mu_k^j$.
The above process starts with some initial beliefs $\mu_{0}^i$, $i=1,\ldots,n$.
Note that, if $\Theta$ is a finite or a countable set, the update rule in Eq. \eqref{protocol_cont} reduces to:
	for every $B\subseteq \Theta$,
	\begin{align}\label{protocol_dis}
	\mu_{k}^i\left(B\right) & = \frac{1}{Z_{k}^i}\sum\limits_{\theta \in B} \prod\limits_{j=1}^{n} \mu_{k-1}^j\left(\theta\right)^{a_{ij}}  \ell^i\left(s_{k}^i|\theta\right)
	\end{align}
The updates  in Eqs.~\eqref{protocol_cont} and~\eqref{protocol_dis} 
can be viewed as two-step processes. First, every agent constructs an aggregate belief using weighted geometric average 
of its own belief and the beliefs of its neighbors. Then, 
each agent performs a Bayes' update using the aggregated belief as a prior.

\subsection{Connection with the Stochastic Mirror Descent Method}
To make this connection\footnote{Particularly simple when $\Theta$ is a finite set.}, we observe that 
the optimization problem in Eq.~\eqref{opt_problem} is equivalent to the following problem:
\begin{align*}
\min_{\pi \in \mathscr{P}_{\Theta} } \mathbb{E}_{\pi}  \sum\limits_{i=1}^n D_{KL}\left(f^i\|\ell^i\right)
& = \min_{\pi \in \mathscr{P}_{\Theta} } \sum\limits_{i=1}^n \mathbb{E}_{\pi} \mathbb{E}_{f^i} [- \log \ell^i]
\end{align*} 
where $\mathscr{P}_{\Theta}$ is the set of all distributions on $\Theta$.
Under some technical conditions the expectations can exchange the order, so the problem in Eq.~\eqref{opt_problem} 
is equivalent to the following one:
\begin{align}\label{main2}
\min_{\pi \in \mathscr{P}_{\Theta} } \sum\limits_{i=1}^n \mathbb{E}_{f^i} \mathbb{E}_{\pi} [- \log \ell^i]
\end{align}
The difficulty in evaluating the objective function in Eq.~\eqref{main2} 
lies in the fact that the distributions $f^i$ are unknown. 
A generic approach to solving such problems is the class of stochastic approximation methods, where the objective is minimized by constructing a sequence of gradient-based iterates where the true gradient of the objective 
(which is not available) is replaced with a gradient sample that is available at the given update time.
A particular method that is relevant here is the stochastic mirror-descent method which would solve the problem 
in Eq.~\eqref{main2}, in a centralized fashion, by constructing a sequence $\{x_k\}$, as follows:
\begin{align}\label{central}
x_{k} & =  \argmin_{y \in X} \left\lbrace \langle g_{k-1}, y\rangle 
+ \frac{1}{\alpha_{k-1}} D_w(y,x_{k-1})\right\rbrace 
\end{align}
where $g_k$ is a noisy realization of the gradient of the objective function 
in Eq.~\eqref{main2} and $D_w(y,x)$ is a Bregman distance function associated with a distance-generating 
function $w$, and $\alpha_k>0$ is the step-size. 
If we take $w(t) = t \log t$ as the distance-generating function, then the corresponding Bregman distance is 
the Kullback-Leiblier (KL) divergence $D_{KL}$. Thus, in this case, 
{\it the update rule in Eq.~\eqref{protocol_cont} corresponds to
a distributed implementation of the stochastic mirror descent algorithm} in~\eqref{central},
where $D_w(y,x)=D_{KL}(y,x)$ and the stepsize is fixed, i.e..,  $\alpha_k =1$ for all $k$.


We summarize the preceding discussion in the following proposition. 

\begin{proposition}\label{mirror}
	The update rule in Eq. \eqref{protocol_cont} {defines a probability measure $\mu_k^i$over the set $\Theta$ generated by the probability density $\bar{\mu}_k^i = d\mu_k^i/d\lambda$  that coincides with the solution of the distributed stochastic mirror descent algorithm applied to the optimization problem in Eq.~\eqref{opt_problem}., i.e.}
	{
		\begin{align}\label{stoc_mirror}
		\bar{\mu}_{k}^i = \argmin_{\pi \in P_{\Theta} } \left\lbrace 
		\mathbb{E}_{\pi} [\textendash\log \ell^i(s_{k}^i|\cdot)] + \sum\limits_{j=1}^{n} a_{ij} D_{KL}(\pi\|\bar{\mu}_{k-1}^j)\right\rbrace 
		\end{align}
	}
\end{proposition}
\begin{IEEEproof}
We need to show that {the density $\bar{\mu}_k^i$ associated with the probability measure $\mu_k^i$} defined by
Eq.~\eqref{protocol_cont} minimizes the problem in Eq.~\eqref{stoc_mirror}.

	First, define the argument in Eq.~\eqref{stoc_mirror} as
	\begin{align*}
	G(\pi) & = - \mathbb{E}_{\pi} \log \ell^i\left(s_{k}^i|\cdot\right) + \sum_{j=1}^{n} a_{ij} D_{KL}\left(\pi||\bar{\mu}_{k-1}^j\right)
	\end{align*}
	and add and subtract the KL divergence between $\pi$ and the density  $\bar{\mu}_{k}^i$ to obtain
	\begin{align*}
	G(\pi) & 
	= - \mathbb{E}_{\pi} \log \ell^i(s_{k}^i|\cdot) + \sum\limits_{j=1}^{n} a_{ij} D_{KL}(\pi\|\bar{\mu}_{k-1}^j)- D_{KL}(\pi\|\mu_{k}^i) + D_{KL}(\pi\|\bar{\mu}_{k}^i) \\
	& =- \mathbb{E}_{\pi} \log \ell^i\left(s_{k}^i|\cdot\right) 
	+ D_{KL}\left(\pi\|\bar{\mu}_{k}^i\right)+ \sum\limits_{j=1}^{n}a_{ij} \mathbb{E}_\pi\left(  \log \frac{\pi}{\bar{\mu}_{k-1}^j} - \log \frac{\pi}{\bar{\mu}_{k}^i}\right)\\
	&= - \mathbb{E}_{\pi} \log \ell^i(s_{k}^i|\cdot) 
	+ D_{KL}\left(\pi\|\bar{\mu}_{k}^i\right) + \sum\limits_{j=1}^{n}a_{ij} \mathbb{E}_\pi \log \frac{\bar{\mu}_{k}^i}{\bar{\mu}_{k-1}^j} 
	\end{align*}
	Now, we use the relation for the density $\bar{\mu}_k^i = d\mu_k^i/d\lambda$,
	which is implied by the update rule for $\mu_{k}^i$ in Eq.~\eqref{protocol_cont}, and obtain
	\begin{align*}
	 G(\pi) 
         &= - \mathbb{E}_{\pi} \log  \ell^i\left(s_{k}^i|\cdot\right)
	+ D_{KL}\left(\pi\|\bar{\mu}_{k}^i\right) + \sum\limits_{j=1}^{n}a_{ij} \mathbb{E}_\pi \log \left(  \frac{1}{\bar{\mu}_{k-1}^j} \frac{1}{Z_{k}^i} \prod\limits_{m=1}^{n} \left(\bar{\mu}_{k-1}^m\right)^{a_{im}}\ell^i(s_{k}^i|\cdot) \right)  \\
	&= -  \log Z_{k}^i + D_{KL}\left(\pi\|\bar{\mu}_{k}^i\right) 
	\end{align*}
	The first term in the preceding line does not depend on the distribution $\pi$. Thus, 
	we conclude that the solution to the problem in Eq.~\eqref{stoc_mirror} is 
	the density $\pi^* = \bar{\mu}_{k}^i$ (almost everywhere).
\end{IEEEproof}

\section{Countable Hypothesis Set}\label{sec_count}
In this section we present a concentration result for 
the update rule in Eq.~\eqref{protocol_dis} specific for the case of a countable hypothesis set. 
Later in Section \ref{sec_cont} we will analyze the case of  $\Theta\subset \mathbb{R}^d$.

First, we provide some definitions that will help us build the desired results. 
Specifically, we will study how the beliefs of all agents concentrate around the true hypothesis~$\theta^*$. 

\begin{definition}\label{kl_balls}
	Define a Kullback-Leibler Ball (KL) of radius $r$ centered at $\theta^*$ as.
	\begin{align*}
	\mathcal{B}_r(\theta^*) & = \left\lbrace \theta \in \Theta \left|   \frac{1}{n} \sum_{i=1}^{n} D_{KL}\left(\ell^i\left(\cdot|\theta^*\right) ,\ell^i\left(\cdot|\theta \right) \right) \right. \leq r  \right\rbrace 
	\end{align*}
\end{definition}

\begin{definition}\label{covering}
	Define the covering of the set  $\mathcal{B}_r^c(\theta^*)$ generated by a strictly increasing sequence $\{r_l\}_{l=1}^{\infty}$ with $r_1=r$ as the union of disjoint KL bands as follows:
	\begin{align*}
	\mathcal{B}_r^c(\theta^*) & = \bigcup_{l=1}^{\infty} \{\mathcal{B}_{r_{l+1}}(\theta^*) \backslash \mathcal{B}_{r_{l}}(\theta^*)\}
	\end{align*}
	where $\{\mathcal{B}_{r_{l+1}} \backslash \mathcal{B}_{r_{l}}\}$ denotes the complement between the set $\mathcal{B}_{r_{l+1}}$ and the set $\mathcal{B}_{r_{l+1}}$, i.e. $\mathcal{B}_{r_{l+1}} \cap \mathcal{B}_{r_{l}}^c$
	We denote the cardinality $\{\mathcal{B}_{r_{l+1}}\backslash \mathcal{B}_{r_{l}} \}$ by $\mathcal{N}_{r_l}$, i.e. $\left| \{\mathcal{B}_{r_{l+1}}\backslash \mathcal{B}_{r_{l}} \} \right| = \mathcal{N}_{r_l}$.
\end{definition}

\begin{figure}[ht]
	\centering
	 \begin{overpic}[width=0.4\textwidth]{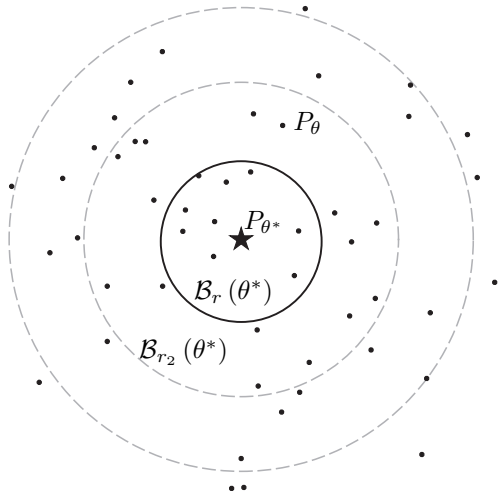}
			\put(38,39){{\small $\mathcal{B}_r\left(\theta^*\right)$}}
			\put(48,53){{\small $P_{\theta^*}$}}
			\put(58,73){{\small $P_{\theta}$}}
			\put(27,27){{\small $\mathcal{B}_{r_2}\left(\theta^*\right)$}}
	 \end{overpic} 
	 \captionsetup{justification=centering}
	 \caption{Creating a covering for a ball $\mathcal{B}_r\left(\theta^*\right)$. $\bigstar$ represents the correct hypothesis 
	 $\ell^i\left(\cdot|\theta^*\right)$, $\bullet$ indicates the location of other hypotheses and the dash lines indicates the boundary of the balls $\mathcal{B}_{r_l}\left(\theta^*\right)$.}
	 \label{lcovering_dis}
\end{figure}

We are interested in bounding the beliefs' concentration on a ball $\mathcal{B}_r\left(\theta^*\right)$ for an arbitrary $r>0$, which is based on a covering of the complement set $\mathcal{B}_r^c\left(\theta^*\right)$.
To this end, Definitions~\ref{kl_balls} and~\ref{covering} provide the tools for constructing such a covering. The strategy is to analyze how the hypotheses are distributed in the space of probability distributions, see Figure \ref{lcovering_dis}. The next assumption will provide conditions on the hypothesis set which guarantee the concentration results.

 \begin{assumption}\label{conv_hyp}
	The series $\Sigma_{l \geq 1}  \exp\left(-r_{l}^2 +\log \mathcal{N}_{r_l}\right)$ converges, 
	where the sequence $\{r_l\}$ is as in Definition~\ref{covering}.
\end{assumption}

We are now ready to state the main result for a countable hypothesis set $\Theta$.

\begin{theorem}\label{main_count}
	Let Assumptions~\ref{assum:init},~\ref{assum:graph} and~\ref{conv_hyp} hold, 
	and let $\rho \in (0,1)$ be a desired probability  tolerance. 
	Then, the belief sequences $\{\mu_{k}^i\}$, $i\in V$, generated by the update rule in Eq.~\eqref{protocol_dis},
	with the initial beliefs such that ${\mu_0^i(\theta^*)>\epsilon}$ for all $i$,  
	have the following property: for any $\sigma \in (0,1)$ and any radius $r>0$ with probability $1-\rho$,
	\begin{align*}
	\mu_{k}^i\left(\mathcal{B}_r\left(\theta^*\right)\right) & \geq 1 - \sigma \qquad \text{ for all } i \text{ and all }k\geq N  
	\end{align*}
	where $N = \min_{k \geq 1} \{k\in N_1\cap N_2\}$ with the sets $N_1$ and $N_2$ given by
	{
	\begin{align*}
	N_1& =\left\lbrace k\mid \exp\left(\frac{ 1}{8\log^2\frac{1}{\alpha}}\right) \sum\limits_{l \geq 1} \mathcal{N}_{r_l} \exp \left(  -kr_{l}^2 \right)  \leq \rho  \right\rbrace \\
	N_2& =\left\lbrace  k\mid { C_3 }\exp\left(-\frac{k}{2}\gamma(\theta)\right) \leq \sigma {\prod_{i=1}^n \mu_0^i\left(\theta\right)^{\frac{1}{n}} },   \forall\theta \not\in \mathcal{B}_r\left(\theta^*\right)  \right\rbrace 
	\end{align*}
}
	where  {	$\gamma(\theta)  =  \frac{1}{n}\sum_{i=1}^{n}D_{KL}\left(\ell^i\left(\cdot|\theta^*\right)||\ell^i\left(\cdot|\theta\right)\right)$, ${C_3  = \frac{1}{\epsilon}\exp\left(\frac{8 \log \frac{1}{\alpha} \log n}{1 - \lambda}\right)}$},
	$\alpha$ is as in Assumption~\ref{assum:init}(b), $\mathcal{N}_{r_l}$ and $r_l$ are as in 
	Definition~\ref{covering}, while 
	${\lambda = 1-\eta/4n^2}$. 
	If $A$ is a lazy-metropolis matrix, then ${\lambda = 1- 1/\mathcal{O}(n^2) }$.
\end{theorem}

Observe that if $k\in N_1$, then $m\in N_1$ for all $m\ge k$, and the same is true for the set $N_2$,
so we can alternatively write
\[N=\max\left\{ \min_{k\ge 1} \{k\in N_1\}, \ \min_{k\ge 1} \{k\in N_2\}\right\}.\]
Further, note that $N$ depends on the radius $r$ of the KL ball, as the set $N_1$ 
involves $\mathcal{N}_{r_{l}}$ and $r_l$ which both depend on~$r$, while 
the set $N_2$ explicitly involves $r$. Finally, note that the smaller the radius $r$, the larger $N$ is. 
We see that $N$ also depends on 
the number $n$ of agents, the learning parameter $\alpha$, 
the learning capabilities of the network represented by $\gamma(\theta)$, the initial beliefs $\mu_0^i$, 
the number of hypotheses that are far away from $\theta^*$ and
their probability distributions.

Theorem~\ref{main_count} states that the beliefs of all agents will concentrate within the KL ball 
$\mathcal{B}_r\left(\theta^*\right)$
with a radius $r>0$ for a large enough $k$, i.e., $k\ge N$.
Note that the (large enough) index $N$ is determined as the smallest $k$ for which two relations are satisfied,
namely, the relations defining the index sets $N_1$ and $N_2$.
The set $N_1$ contains all indices $k$ for which a weighted sum of the total mass of the hypotheses 
$\theta\notin
\mathcal{B}_r\left(\theta^*\right)$ is small enough (smaller than the desired probability tolerance $\rho$). 
Specifically, we require the number $\mathcal{N}_{r_l}$ of hypothesis in the $l$-th band does not grow 
faster than the squared radius $r_{l}^2$ of the band, i.e., the wrong hypothesis should not accumulate too fast far away from the true hypothesis $\theta^*$. Moreover, the condition in $N_1$ also prevents having an infinite number of hypothesis per band. The set $N_2$ captures the iterations $k$ at which, for all agents,  the current beliefs $\mu_k^i$ had recovered 
from the the cumulative effect of ``wrong" initial beliefs that had given probability masses 
to hypotheses far away from $\theta^*$.

In the proof for Theorem \ref{main_count}, we use the relation between
the posterior beliefs and the initial beliefs on a measurable set $B$ such that $\theta^* \in B$. For such a set, we have
{
\begin{align}\label{initial_measure}
\mu_{k}^i\left(B\right) & = \frac{1}{Z_{k}^i}\sum\limits_{\theta \in B} \prod\limits_{j=1}^{n} \mu_0^j\left(\theta\right)^{\left[A^{k}\right]_{ij} } \prod\limits_{t=1}^{k} \prod\limits_{j=1}^{n}  \ell^j(s_{t}^j|\theta)^{\left[A^{k-t}\right]_{ij}}
\end{align}
}
where $ Z_{k}^i $ is the appropriate normalization constant. 
Furthermore, after a few algebraic operations we obtain
{
\begin{align}\label{ready_for_bounds}
\mu_{k}^i\left(B\right) & \geq 1\text{--}\sum\limits_{\theta \in B^c} { \prod\limits_{j=1}^{n} \left( \frac{\mu_0^j(\theta)}{\mu_0^j(\theta^*)}\right) ^{\left[A^{k}\right]_{ij} }} \prod\limits_{t=1}^{k} \prod\limits_{j=1}^{n} \left(  \frac{\ell^j(s_{t}^j|\theta)}{\ell^j(s_{t}^j|\theta^*)}\right) ^{\left[A^{k\text{-}t}\right]_{ij}}
\end{align}
}
{
Moreover, since $\mu_0^j(\theta^*) > \epsilon$ for all $j$, it follows that}
{
	\begin{align}\label{ready_for_bounds2}
	\mu_{k}^i\left(B\right) & \geq 1\text{--}\frac{1}{\epsilon}\sum\limits_{\theta \in B^c}\prod\limits_{t=1}^{k} \prod\limits_{j=1}^{n} \left(  \frac{\ell^j(s_{t}^j|\theta)}{\ell^j(s_{t}^j|\theta^*)}\right) ^{\left[A^{k\text{-}t}\right]_{ij}}
	\end{align}
}

Now we will state a useful result from \cite{sha14} which will allow us to bound the right hand term of Eq. \eqref{ready_for_bounds}.
\begin{lemma}\label{shahin}[Lemma 2 in \cite{sha14}]
	Let Assumptions \ref{assum:graph} hold, then the matrix $A$ satisfies: for all $i$, 
	\begin{align*}
	\sum\limits_{t=1}^{k} \sum\limits_{j=1}^{n} \left| \left[A^{k-t}\right]_{ij} - \frac{1}{n} \right| & \leq \frac{4 \log n}{1-\lambda} 
	\end{align*}
	where $\lambda = 1-\eta / 4n^2$, and if $A$ is a lazy-metropolis matrix associated with $\mathcal{G}$ then $\lambda = 1- 1 / \mathcal{O}(n^2)$.
\end{lemma}
{
	If follows from Eq. \eqref{ready_for_bounds2}, Lemma \ref{shahin} and Assumption \ref{assum:init} that
	{
		\begin{align}\label{bounded2}
		\mu_{k}^i\left(B\right) & \geq 1\textendash  C_3\sum\limits_{\theta \in B^c} \prod\limits_{t=1}^{k} \prod\limits_{j=1}^{n} \left(  \frac{\ell^j(s_{t}^j|\theta)}{\ell^j(s_{t}^j|\theta^*)}\right) ^{\frac{1}{n}}
		\end{align}
	}
	for all $i$, where $C_3$ is as defined in Theorem \ref{main_count}.
 }
Next we provide an auxiliary result about the concentration properties of the beliefs on a set $B$.
\begin{lemma}\label{bound_prob}
	For any $k\ge0$ it holds that
	\begin{align*}
\mathbb{P}_{\boldsymbol{f}}\left(\bigcup_{\theta \in B^c} \left\lbrace\bar{v}_{k}(\theta)   \geq -\frac{k}{2} \gamma(\theta) \right\rbrace  \right) 
& \leq C_2\sum\limits_{l \geq 1} \mathcal{N}_{r_l} \exp \left(  -kr_{l}^2 \right)
	\end{align*}
	where $C_2 = \exp\left(\frac{1}{8} \frac{ 1}{\log^2\frac{1}{\alpha}}\right) $, $\gamma(\theta)$ and $\mathcal{N}_{r_l}$ and $r_l$ are as in Theorem \ref{main_count}.
\end{lemma}

\begin{IEEEproof}
	First define the following random variable
	\begin{align*}
	\bar{v}_{k}(\theta) 
	& = \sum\limits_{{t=1}}^{k} \frac{1}{n}\sum\limits_{i=1}^{n} \log \frac{\ell^i(S_{t}^i|\theta)}{\ell^i(S_{t}^i|\theta^*)}  
	\end{align*}
	Then, by using the union bound and McDiarmid inequality we have,
	{
\begin{align*}
 \mathbb{P}_{\boldsymbol{f}}\left(\bigcup_{\theta \in B^c} \left\lbrace\bar{v}_{k}(\theta) 
 \text{-} \mathbb{E}_{\boldsymbol{f}}\bar{v}_{k}(\theta)   
 \geq \bar{\epsilon} \right\rbrace  \right)
 &  \leq \sum\limits_{\theta \in B^c} \mathbb{P}_{\boldsymbol{f}}
 \left\{ \bar{v}_{k}(\theta) - \mathbb{E}_{\boldsymbol{f}}\bar{v}_{k}(\theta) \geq  \bar{\epsilon} \right\}\\
 & \leq \sum\limits_{\theta \in B^c} \exp \left(  \frac{2\bar{\epsilon}^2}{ 4k \log^2\frac{1}{\alpha} }\right)
\end{align*}
}
and by setting $\epsilon =-\frac{1}{2}\mathbb{E}_{\boldsymbol{f}}\left[\bar{v}_{k}(\theta)\right]$, it follows that
{
\begin{align*}
 \mathbb{P}_{\boldsymbol{f}}
 \left(\bigcup_{\theta \in B^c} \left\lbrace \bar{v}_{k}(\theta)   
 \geq \frac{1}{2} \mathbb{E}_{\boldsymbol{f}}\left[\bar{v}_{k}(\theta) \right]\right\rbrace  \right)
  \leq 
 \sum\limits_{\theta \in B^c} \exp \left( \frac{\left( \mathbb{E}_{\boldsymbol{f}}\left[\bar{v}_{k}(\theta)\right]\right) ^2}{8k\log^2\frac{1}{\alpha}}  \right)
\end{align*}
}
It can be seen that  $\mathbb{E}_{\boldsymbol{f}}\left[\bar{v}_{k}(\theta)\right] = - k\gamma(\theta)$, thus yielding
\begin{align*}
\mathbb{P}_{\boldsymbol{f}}
\left(\bigcup_{\theta \in B^c} \left\lbrace\bar{v}_{k}(\theta)   \geq -\frac{k}{2} \gamma(\theta) \right\rbrace  \right) 
& \leq C_2\sum\limits_{\theta \in B^c}  \exp \left(  -k \gamma^2(\theta)  \right)
\end{align*}
Now, we let the set $B$ be the KL ball of a radius $r$ centered at $\theta^*$ 
and follow Definition~\ref{covering} to exploit the representation of 
$\mathcal{B}_r^c(\theta^*)$ as the union of KL bands, for which we obtain
{
\begin{align*}
\sum_{\theta \in \mathcal{B}_{r}^c(\theta^*)} \exp \left(  -k\gamma^2(\theta) \right) 
&= \sum\limits_{l \geq 1} \sum_{\theta \in \mathcal{B}_{r_{l+1}} \backslash \mathcal{B}_{r_{l}}} \exp \left(  -k\gamma^2(\theta) \right)\\
&  \leq \sum\limits_{l \geq 1} \mathcal{N}_{r_l} \exp \left(  -kr_{l}^2 \right)
\end{align*}
}
thus, completing the proof.
\end{IEEEproof}
We are now ready to proof Theorem \ref{main_count}

\begin{IEEEproof}[Theorem \ref{main_count}]
	From Lemma~\ref{bound_prob}, 
	where we take $k$ large enough to ensure the desired probability, it follows that
	with probability $1-\rho$, we have: for all $k\in N_1$,
	\begin{align*}
	\mu_{k}^i\left(\mathcal{B}_r(\theta^*)\right) & 
	\geq 1-C_3\sum\limits_{\theta \in \mathcal{B}_r^c(\theta^*) }\exp\left(-\frac{k}{2}\gamma(\theta)\right)\\
	& \geq 1-\sum\limits_{\theta \in \mathcal{B}_r^c(\theta^*) } \sigma {\prod_{i=1}^n \mu_0^i\left(\theta\right)^{\frac{1}{n}} } \\
	& \geq 1-\sigma
	\end{align*}
	where the last inequality follows from Eq.~\eqref{bounded2} where we further take sufficiently large $k$.
\end{IEEEproof}

\section{Continuum of Hypotheses}\label{sec_cont}
In this section we will provide the concentration results 
for a continuous hypothesis set $\Theta \subseteq \mathbb{R}^d$. 
At first, we present some definitions that we use in constructing coverings analogously to that in Section~\ref{sec_count}. 
In this case, however, we employ the Hellinger distance.
\begin{definition}\label{h_balls}
	Define a Hellinger Ball (H) of radius $r$ centered at $\theta^*$ as.
	\begin{align*}
	\mathcal{B}_r(\theta^*) 
	& = \left\lbrace \theta \in \Theta \left|   \frac{1}{\sqrt{n}}  h\left(\boldsymbol{\ell}\left(\cdot|\theta^*\right) ,\boldsymbol{\ell}\left(\cdot|\theta \right) \right) \right. \leq r  \right\rbrace 
	\end{align*}
\end{definition}
\begin{definition}
	Let $\left(M,d\right)$ be a metric space. A subset $S_{\delta} \subseteq M$ is called $\delta$-separated with $\delta>0$ if $d(x,y)\geq \delta$ for any $x,y \in S_{\delta}$. Moreover, for a set 
	$B \subseteq M$, let $\mathcal{N}_B(\delta)$ be the smallest number of Hellinger balls with centers in $S_\delta$ 
	of radius $\delta>0$ needed to cover the set $B$, i.e., such that $B \subseteq \bigcup_{ z \in S_\delta} \mathcal{B}_{\delta}\left(z\right)$.
\end{definition}

\begin{definition}\label{covering_cont}
Let $\{r_l\}$ be a strictly decreasing sequence such that ${r_1 =1}$ and $\lim_{l \to \infty} r_l = 0$.
	Define the covering of the set $\mathcal{B}_{r}^c(\theta^*)$ generated by the sequence $\{r_l\}$ as follows:
	\begin{align*}
	\mathcal{B}_{r}^c(\theta^*)& = \bigcup_{l = 1}^{{{L_r}}-1} \{\mathcal{B}_{r_{l}} \backslash\mathcal{B}_{r_{l+1}} \}
	\end{align*}
	where ${L_r}$ is the smallest $l$ such that $r_l \leq r$.
	Moreover, given a positive sequence $\{\delta_l\}$, 
	we denote by $\mathcal{N}_{r_{l}}(\delta_l)$ the maximal $\delta_l$-separated subset of the set $\{\mathcal{B}_{r_l} \backslash\mathcal{B}_{r_{l+1}} \}$ and denote its cardinality by 
	$K_l$, i.e. $K_l = |\mathcal{N}_l(\delta_l)|$. 
	Therefore, we have the following covering of $\mathcal{B}_{r}^c(\theta^*)$,
	\begin{align*}
	\mathcal{B}_{r}^c(\theta^*) & = \bigcup_{l \ge 1}^{{{L_r}}-1} \bigcup_{z_m\in \mathcal{N}_l(\delta_l)} \mathcal{F}_{l,m} 
	\end{align*}
	where $\mathcal{F}_{l,m} = \mathcal{B}_{\delta_l}(z_m \in \mathcal{N}_{r_{l}}(\delta_l)) \cap \{\mathcal{B}_{r_{l}} \backslash\mathcal{B}_{r_{l+1}} \} $.
\end{definition}

\begin{figure}[ht]
	\centering
	\begin{overpic}[width=0.4\textwidth]{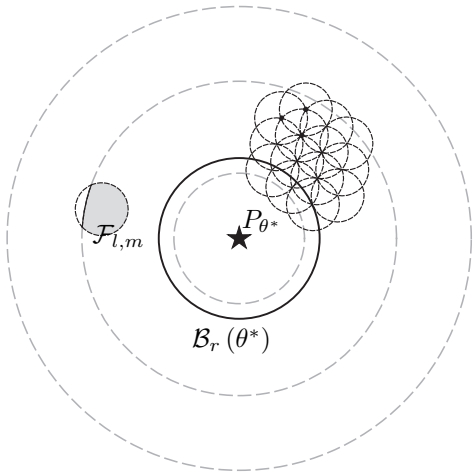}
		\put(38,29){{\small $\mathcal{B}_{r}\left(\theta^*\right)$}}
		\put(18,50){{\small $\mathcal{F}_{l,m}$}}
		\put(48,53){{\small $P_{\theta^*}$}}
	\end{overpic} 
	\captionsetup{justification=centering}
	\caption{Creating a covering for a set $\mathcal{B}_r\left(\theta^*\right)$. $\bigstar$ represents the correct hypothesis $\ell^i\left(\cdot|\theta^*\right)$. }
	\label{lcovering_cont}
\end{figure}

Figure \ref{lcovering_cont} depicts the elements of a covering for a set $\mathcal{B}_r^c\left(\theta^*\right)$. The cluster of circles at the top right corner represents the balls $ \mathcal{B}_{\delta_l}(z_m \in \mathcal{N}_{r_{l}}(\delta_l))$ and for a specific case in the left of the image we illustrate the set $\mathcal{F}_{l,m}$.

%
We are now ready to state the main result regarding continuous set of hypotheses $\Theta \subseteq \mathbb{R}^d$.

\begin{theorem}\label{main_cont}
	Let Assumptions~\ref{assum:init},~\ref{assum:graph},  and~\ref{conv_hyp} hold, and let ${\rho \in (0,1)}$ be a given probability tolerance level. 
	Then, the beliefs $\{\mu_k^i\},$ $i\in V,$ generated by the update rule in~Eq.~\eqref{protocol_cont} 
	with uniform initial beliefs, 
	are such that, for any $\sigma \in (0,1)$ and any $r>0$ with probability $1-\rho$,
	\begin{align*}
	\mu_{k}^i\left(\mathcal{B}_{r}\left(\theta^*\right)\right) & \geq 1 - \sigma \qquad \text{ for all } i \text{ and all }k\geq N  
	\end{align*}
	where $N = \min \left\lbrace k \geq 1 \left|  k \in N_1 \text{ and } k \in N_2 \right. \right\rbrace $ with
	{
	\begin{align*}
	N_1 &\text{=}\left\lbrace  \sum_{l=1}^{{L_r}-1}  \exp\left(\frac{8 \log \frac{1}{\alpha} \log n}{1 - \lambda} \text{ --}k \left( r_{l+1}\text{ --}\delta_l \text{ --}R\right) \text{--}d \log \delta_l \right)\leq \rho  \right\rbrace \\
	N_2 & \text{=}\left\lbrace  \sum_{l=1}^{{L_r}-1} \exp \left(d\log\frac{r_l}{R} - 2k\left(r_{l+1}-\delta_l -R\right)\right) \leq \sigma  \right\rbrace 
	\end{align*}
	}
	for a parameter $R$ such that $r>R$ and  $r_{l+1}-\delta_l -R >0$ for all $l\ge1$. 
	The constant $\alpha$ is as in Assumption~\ref{assum:init}, $d$ is the dimension of the space of $\Theta$, 
	$r_{l}$ and $\delta_l$ are as in Definition~\ref{covering_cont}, while $\lambda$ is the same as in 
	Theorem~\ref{main_count}. 
\end{theorem}
Analogous to Theorem \ref{main_count}, 
Theorem \ref{main_cont} provides a probabilistic concentration result for the agents' beliefs around 
a Hellinger ball of radius $r$ with center at $\theta^*$ for sufficiently large $k$. 

Similarly to the preceding section, we represent the beliefs $\mu_{k}^i$  in terms of the initial beliefs 
and the cumulative product of the weighted likelihoods received from the neighbors. 
In particular, analogous to Eq.~\eqref{initial_measure}, we have that for every $i$ and 
for every measurable set $B \subseteq \Theta$:
\begin{align}\label{initial_measure_c}
\mu_{k}^i\left(B\right) & 
= \frac{1}{Z_{k}^i}\int\limits_{\theta \in B} \prod\limits_{t=1}^{k} \prod\limits_{j=1}^{n} 
\ell^j(s_{t}^j|\theta)^{\left[A^{k-t}\right]_{ij}}    d\mu_0^j\left(\theta\right)  
\end{align}
with the corresponding normalization constant $Z_k^i$, and assuming all agents have uniform beliefs at time $0$.

It will be easier to work with the beliefs' densities, so we define the density of a measurable set with respect to the observed data. 
\begin{definition}\label{g_densities}
	The density $g_B^i$ of a measurable set $B\subseteq \Theta$, where $\mu_0^i\left(B\right) > 0$ with respect to the product distribution of the observed data is given by
	\begin{align}\label{def_dens}
	g_B^i\left(\hat{s}\right) & = \frac{1}{\mu_0^i\left(B\right)} 
	    \int\limits_{B} \prod\limits_{t=1}^{k} 
	\prod\limits_{j=1}^{n} \ell^j(S_{t}^j|\theta)^{\left[A^{t-k}\right]_{ij}}d\mu_0^j\left(\theta\right)
	\end{align}
	where $\hat{s} = \{S_t^i\}_{t=1:k}^{i=1:n}$ and $P_B^i = g_B^i \cdot (\lambda^i)^{\otimes k}$. 
\end{definition}

The next lemma relates the density $g_B^i\left(\hat{s}\right)$ which is defined per agent to a quantity that is common among all nodes in the network.
\begin{lemma}\label{dist_densities}
	Consider the densities as defined in Eq. \eqref{def_dens}, then
	\begin{align}
	g_B^i\left(\hat{s}\right) & \leq C_1g_B\left(\hat{s}\right) ^{\frac{1}{n}}
	\end{align}
	where $C_1 = \exp\left(\frac{8 \log \frac{1}{\alpha} \log n}{1 - \lambda}\right)   $ and
	\begin{align*}
	g_B\left(\hat{s}\right)
	 & = \int\limits_{B}  \prod\limits_{t=1}^{k}  \prod\limits_{j=1}^{n}  \ell^j(S_{t}^j|\theta)d\mu_0^j\left(\theta\right).
	\end{align*}
\end{lemma}
\begin{IEEEproof}
	By definition of the densities, we have
	{
	\begin{align*}
	 g_B^i\left(\hat{s}\right) & = \frac{1}{\mu_0^i\left(B\right)} \int\limits_{B} \prod\limits_{t=1}^{k} \prod\limits_{j=1}^{n} \ell^j(S_{t}^j|\theta)^{\left[A^{t-k}\right]_{ij}}d\mu_0^j\left(\theta\right)\\
	& = \frac{1}{\mu_0^i\left(B\right)} \int\limits_{B} \exp\left(\sum\limits_{{t=1}}^{k}\sum\limits_{j=1}^{n} \left[A^{t-k}\right]_{ij} \log \ell^j(S_{t}^j|\theta)\right)d\mu_0^j\left(\theta\right)\\
	& = \frac{1}{\mu_0^i\left(B\right)} \int\limits_{B} \exp\left(\sum\limits_{{t=1}}^{k}\sum\limits_{j=1}^{n} \left( \left[A^{t-k}\right]_{ij} - \frac{1}{n}\right)  \log \ell^j(S_{t}^i|\theta)
	+ \sum\limits_{{t=1}}^{k}\frac{1}{n}  \sum\limits_{j=1}^{n} \log \ell^j(S_{t}^j|\theta) \right)d\mu_0^j\left(\theta\right)
	\end{align*}
	}
	where the last line follows by adding and subtracting $1/n$. Hence, by Lemma~\ref{shahin}, we further obtain
	{
	\begin{align*}
	g_B^i\left(\hat{s}\right) 
	& \leq \frac{C_1}{\mu_0^i\left(B\right)}   \int\limits_{B} \exp\left( \sum\limits_{t=1}^{k}\frac{1}{n}  \sum\limits_{j=1}^{n} \log \ell^j(S_{t}^j|\theta) \right)d\mu_0^j\left(\theta\right) \\ 
	& =  \frac{C_1}{\mu_0^i\left(B\right)}  \int\limits_{B} \prod\limits_{t=1}^{k}  \prod\limits_{j=1}^{n}  \ell^j(S_{t}^j|\theta)^{1/n} d\mu_0^j\left(\theta\right)\\
	& \leq\frac{C_1}{\mu_0^i\left(B\right)} \left(\int\limits_{B}  \prod\limits_{t=1}^{k}  \prod\limits_{j=1}^{n}  \ell^j(S_{t}^j|\theta)d\mu_0^j\left(\theta\right)\right) ^{1/n} 
	\end{align*}
	}
	where the last inequality follows from Jensen's inequality.
\end{IEEEproof}

The next Lemma is an analog of Lemma \ref{bound_prob} which we use to bound the probability concentrations 
with respect to the ratio $g_{\mathcal{F}_{l,m}}^i\left(\hat{s}\right)/ g_{\mathcal{B}_{R}(\theta^*)}^i\left(\hat{s}\right)$.
\begin{lemma}\label{bound_cont}
Consider the ratio $g_{\mathcal{F}_{l,m}}^i\left(\hat{s}\right)/ g_{\mathcal{B}_{R}(\theta^*)}^i\left(\hat{s}\right)$, then
{
\begin{align*}
\mathbb{P}_{\mathcal{B}_{R}(\theta^*)}\left( \log\left( \frac{g_{\mathcal{F}_{l,m}}^i\left(\hat{s}\right)}{g_{\mathcal{B}_{R}(\theta^*)}^i\left(\hat{s}\right)}\right)\geq -2k \left( r_l-\delta_l -R\right) \right) \leq  C_2 \exp\left( -k \left( r_{l+1}-\delta_l -R\right) +d \log \delta_l \right).
\end{align*}
}
with $C_2$ as defined in Lemma \ref{bound_prob}.
\end{lemma}
\begin{IEEEproof}
	By using the union bound, the Markov inequality and Lemma 1 in \cite{bir15}, we have that
	{
	\begin{align*}
	\mathbb{P}_{\mathcal{B}_{R}(\theta^*)}\left(\log\left( \frac{g_{\mathcal{F}_{l,m}}^i\left(\hat{s}\right)}{g_{\mathcal{B}_{R}(\theta^*)}^i\left(\hat{s}\right)}\right)\geq y \right) 
& \leq C_1\exp \left(-\frac{y}{2}\right)  \exp \left( -k\frac{1}{n}h^2\left(\mathcal{F}_{l,m},\mathcal{B}_{R}(\theta^*)\right)\right) \\
	& \leq C_1   \exp\left(-\frac{y}{2} - 2k\left( r_{l+1}-\delta_l -R\right) \right)
	\end{align*}
}	
	where the last inequality follows from Proposition 5 and Corollary 1 in \cite{bir15}, where
	\begin{align*}
	\frac{1}{\sqrt{n}} h\left(\mathcal{F}_{l,m},\mathcal{B}_{R}(\theta^*)\right) & \geq \frac{1}{\sqrt{n}} h\left(z_m \in \mathcal{N}_{r_l}\left(\delta_l\right),\theta^*\right) -\delta_l - R\\
	& \geq  \left(r_{l+1}-\delta_l -R\right) 
	\end{align*}
	The desired result is obtained by letting ${y = -2k \left(r_{l+1}-\delta_l -R\right)}$.
\end{IEEEproof}
	Furthermore,
	{
	\begin{align*}
	\mathbb{P}_{\mathcal{B}_{R}(\theta^*)}\left(\bigcup_{\mathcal{F}_{l,m}}\left\lbrace \log\left( \frac{g_{\mathcal{F}_{l,m}}^i\left(\hat{s}\right)}{g_{\mathcal{B}_{R}(\theta^*)}^i\left(\hat{s}\right)}\right)\geq y\right\rbrace \right) 
	& \leq \sum\limits_{l=1}^{{{L_r}}-1} \sum\limits_{m=1}^{K_l}  \mathbb{P}_{\mathcal{B}_{R}(\theta^*)}\left( \log\left( \frac{g_{\mathcal{F}_{l,m}}^i\left(\hat{s}\right)}{g_{\mathcal{B}_{R}(\theta^*)}^i\left(\hat{s}\right)}\right)\geq y\right) \\
	& \leq   C_2\sum\limits_{l=1}^{{{L_r}}-1} K_l  \exp\left(-k\left(r_{l+1} -\delta_l -R\right) \right)\\
	& \leq  C_2 \sum\limits_{l=1}^{{{L_r}}-1}  \exp \left(-k\left(r_{l+1}-\delta_l -R\right) -d \log \delta_l\right)
	\end{align*}
}
	where the last inequality follows from ${K_l \geq \delta_l^{-d}}$, see~\cite{dum07}, \cite{rog57}. 

Lemma \ref{dist_densities} allows us to represent the beliefs on the set $B^c$ as the cumulative beliefs with respect to the density $g_B\left(\hat{s}\right)$. For this, similarly as in Section \ref{sec_count} we will partition the set $B^c$ into Hellinger bands. Then for each band we will find a covering of $\delta$-separated balls and compute the concentration of probability measure with respect to density rations. 

\begin{IEEEproof}\label{proof_theo_cont}[Theorem \ref{main_cont}]
	Lets consider the Hellinger ball $\mathcal{B}_{r}(\theta^*)$. We thus  have
	{
	\begin{align}\label{mu_density}
		\mu_{k}^i\left(\mathcal{B}_{r}(\theta^*)\right) & = \frac{\int\limits_{\theta \in \mathcal{B}_{r}(\theta^*)}   \prod\limits_{t=1}^{k} \prod\limits_{j=1}^{n} 
			\ell^j(s_{t}^j|\theta)^{\left[A^{k-t}\right]_{ij}}    d\mu_0^j\left(\theta\right)}{\int\limits_{\theta \in\Theta} \prod\limits_{t=1}^{k} \prod\limits_{j=1}^{n} 
			\ell^j(s_{t}^j|\theta)^{\left[A^{k-t}\right]_{ij}}    d\mu_0^j\left(\theta\right)} \nonumber \\
		& \geq 1 - \frac{ \int\limits_{\theta \in \mathcal{B}_{r}^c(\theta^*)} \prod\limits_{t=1}^{k} \prod\limits_{j=1}^{n} 
			\ell^j(s_{t}^j|\theta)^{\left[A^{k-t}\right]_{ij}}    d\mu_0^j\left(\theta\right)}{\int\limits_{\theta \in \mathcal{B}_{R}(\theta^*)}  \prod\limits_{t=1}^{k} \prod\limits_{j=1}^{n} 
			\ell^j(s_{t}^j|\theta)^{\left[A^{k-t}\right]_{ij}}   d\mu_0^j\left(\theta\right)} \nonumber \\
		& \geq 1 - \frac{ \int\limits_{\theta \in \mathcal{B}_{r}^c(\theta^*)}\prod\limits_{t=1}^{k} \prod\limits_{j=1}^{n} 
			\ell^j(s_{t}^j|\theta)^{\left[A^{k-t}\right]_{ij}}    d\mu_0^j\left(\theta\right) }{\mu_0^i\left(\mathcal{B}_{R}(\theta^*)\right) g_{\mathcal{B}_{R}(\theta^*)}^i\left(\hat{s}\right) }
	\end{align}
	}
The construction of a partition of the set $\mathcal{B}_{r}^c$ presented in Definition \ref{covering_cont} 
allows us to rewrite Eq. \eqref{mu_density} as follows: 
{
\begin{align*}
\mu_{k}^i\left(\mathcal{B}_{r}^c(\theta^*)\right) & \leq \frac{\sum\limits_{l=1}^{{{L_r}}-1}\sum\limits_{m=1}^{K_l} \int\limits_{\mathcal{F}_{l,m}} \prod\limits_{t=1}^{k} \prod\limits_{j=1}^{n} 
	\ell^j(s_{t}^j|\theta)^{\left[A^{k-t}\right]_{ij}}d\mu_0^j\left(\theta\right)  }{\mu_0^i\left(\mathcal{B}_{R}(\theta^*)\right) g_{\mathcal{B}_{R}(\theta^*)}^i\left(\hat{s}\right) } \\
& =  \sum\limits_{l=1}^{{{L_r}}-1}\sum\limits_{m=1}^{K_l} \frac{\mu_0^j\left(\mathcal{F}_{l,m}\right) }{\mu_0^i\left(\mathcal{B}_{R}(\theta^*)\right)}  \frac{g_{\mathcal{F}_{l,m}}^i\left(\hat{s}\right)}{g_{\mathcal{B}_{R}(\theta^*)}^i\left(\hat{s}\right)}
\end{align*}
}
Finally by applying Lemma \ref{bound_cont} and the fact that of all agents have uniform initial beliefs  
\begin{align*}
\mu_{k}^i\left(\mathcal{B}_{r}^c(\theta^*)\right)  &  \leq  \sum\limits_{l=1}^{{{L_r}}-1}\sum\limits_{m=1}^{K_l}\frac{ \mu_0^j\left(\mathcal{F}_{l,m}\right)  }{\mu_0^j\left(\mathcal{B}_{R}(\theta^*)\right)}\exp(- 2k\left(r_{l+1}\text{-}\delta_l\text{-}R\right))\\
& \leq  \sum\limits_{l=1}^{{{L_r}}-1}\frac{\mu_0^j\left(\mathcal{B}_{r_{l+1}}(\theta^*)\right)}{\mu_0^j\left(\mathcal{B}_{R}(\theta^*)\right)}  \exp(- 2k\left(r_{l+1}-\delta_l -R\right)) \\
& \leq \sum\limits_{l=1}^{{{L_r}}-1} \exp \left(d \log\frac{r_{l+1}}{R} - 2k\left(r_{l+1}-\delta_l -R\right)\right)
\end{align*}
The last 
inequality follows from the initial beliefs being uniform and the volume ratio of the two Hellinger balls with radius $r_l$ and $R$.
\end{IEEEproof}
\section{Conclusions}\label{sec_conclusion}

We proposed an algorithm for distributed learning with a countable and a continuous sets of hypotheses. 
Our results show non-asymptotic  geometric convergence rates for the concentration of the beliefs around the true hypothesis. 

While the proposed algorithm is motivated by the non-Bayesian learning models, 
we have shown that it is also a specific instance of a distributed stochastic mirror descent applied to a well defined optimization problem consisting of the minimization of the sum of Kullback-Liebler divergences.
This indicates an interesting connection between two ``separate" streams of literature and 
provides an initial step to the study of distributed algorithms in a more general form. 
Specifically, it is interesting to explore how variations on stochastic approximation algorithms will induce new non-Bayesian update rules for more general problems. In particular, 
one would be interested in acceleration results for proximal methods, other Bregman distances and 
other constraints in the space of probability distributions. 

Interaction between the agents is modeled as exchange of local probability distributions (\textit{beliefs}) over the hypothesis set between connected nodes in a graph. This will in general generate high communication loads. Nevertheless, results are an initial study towards the distributed learning problems for general hypothesis sets. 
Future work will consider the effect of parametric approximation of the beliefs such that one only needs to communicate a finite number of parameters such as, for example, in Gaussian Mixture Models or Particle Filters.
\bibliographystyle{IEEEtran} 

\bibliography{IEEEfull,bayes_cons_2}

\end{document}